\documentclass[12pt]{article}
\usepackage{e-jc}
\usepackage{amsmath}
\usepackage{amssymb}
\usepackage{graphicx}
\usepackage{float}
\usepackage{ulem}
\usepackage{color}

\newtheorem{theorem}{Theorem}

\newcommand{\Z}{{\mathbb Z}}

\newtheorem{lemma}[theorem]{Lemma}

\newtheorem{claim}{Claim}
\newtheorem{fact}{Fact}

\newcommand{\textdef}{\textit}

\newenvironment{proof}{\noindent{\bf Proof.\,}}{\hfill$\Box$}
\newenvironment{proofcite}[1]{\noindent{\bf Proof of #1.\,}}{\hfill$\Box$}

\begin{document}

\title{Lower bounds for identifying codes in some infinite grids}
\author{Ryan Martin\thanks{Research supported in part by NSA grant H98230-08-1-0015 and NSF grant DMS 0901008}\qquad Brendon Stanton\\Department of Mathematics\\Iowa State University\\Ames, IA 50010}
\maketitle

\begin{abstract}An $r$-identifying code on a graph $G$ is a set $C\subset V(G)$ such that for every vertex in $V(G)$, the intersection of the radius-$r$ closed neighborhood with $C$ is nonempty and unique.  On a finite graph, the density of a code is $|C|/|V(G)|$, which naturally extends to a definition of density in certain infinite graphs which are locally finite.  We present new lower bounds for densities of codes for some small values of $r$ in both the square and hexagonal grids.
\end{abstract}

\section{Introduction}

Given a connected, undirected graph $G=(V,E)$, we define
$B_r(v)$--called the ball of radius $r$ centered at $v$ to be
$$B_r(v)=\{u\in V(G): d(u,v)\le r\}. $$

A subset $C$ of $V(G)$ is called an \textdef{$r$-identifying code} (or simply a \textdef{code}, when $r$ is understood) if it has the properties:
\begin{enumerate}
\item $B_r(v) \cap C \neq \emptyset$
\item $B_r(u) \cap C \neq B_r(v)\cap C$, for all $u\neq v$
\end{enumerate}
The elements of a code $C$ are called \textdef{codewords}.
When $C$ is understood, we define $I_r(v)=I_r(v,C)=B_r(v)\cap C$.  We call
$I_r(v)$ the identifying set of $v$.

Vertex identifying codes were introduced in~\cite{Karpovsky1998} as a way to help with fault diagnosis in multiprocessor computer
systems.  Codes have been studied in many graphs, but of particular interest are codes in the infinite triangular, square, and hexagonal lattices as well as the square lattice with diagonals (king grid). For each of these graphs, there is a characterization so that the vertex set is $\Z\times\Z$.  Let $Q_m$ denote the set of vertices $(x,y)\subset \Z\times\Z$ with $|x|\le m$ and $|y|\le m$.  We may then define the density of a code $C$ by
$$D(C)=\limsup_{m\rightarrow\infty}\frac{|C\cap Q_m|}{|Q_m|}.$$

Our first two theorems, Theorem~\ref{hexr2theorem} and Theorem~\ref{squarer2theorem}, rely on a key lemma, Lemma~\ref{generalpairlemma}, which gives a lower bound for the density of a code assuming that we are able to show that no codeword appears in ``too many'' identifying sets of size 2.  Theorem~\ref{hexr2theorem} follows immediately from Lemma~\ref{generalpairlemma} and Lemma~\ref{hex2lemma} while Theorem~\ref{squarer2theorem} follows immediately from Lemma~\ref{generalpairlemma} and Lemma~\ref{squarepairlemma1}.

    \begin{theorem}\label{hexr2theorem} The minimum density of a 2-identifying code of the
    hex grid is at least 1/5.
    \end{theorem}

    \begin{theorem}\label{squarer2theorem}
    The minimum density of a 2-identifying code of the
    square grid is at least $3/19\approx 0.1579$.
    \end{theorem}

Theorem~\ref{squarer2theorem} can be improved via Lemma~\ref{strongsqpairlemma}, which has a more detailed and technical proof than the prior lemmas.  The idea the lemma is that even though it is possible for a codeword to be in 8 identifying sets of size 2, this forces other potentially undesirable things to happen in the code.  We use the discharging method to show that on average a codeword can be involved in no more than 7 identifying sets of size 2. Lemma~\ref{strongsqpairlemma} leads to the improvement given in Theorem~\ref{squarer2theorem}.

   \begin{theorem}\label{theorem:mainsquare}
     The minimum density of a 2-identifying code of the
     square grid is at least $6/37\approx 0.1622$.
   \end{theorem}

The paper is organized as follows: Section~\ref{generalsection} focuses on some key definitions that we use throughout the paper, provides the proof of Lemma~\ref{generalpairlemma} and provides some other basic facts.  Section~\ref{hexsection} states and proves Lemma~\ref{hex2lemma} from which Theorem~\ref{hexr2theorem} immediately follows.  It is possible to also use this technique to show that the density of a 3-identifying code is at least 3/25, but the proof is long and the improvement is minor so we will exclude it here.  (The proof of this fact will appear in the second author's dissertation~\cite{StantonPending}).  Section~\ref{squaresection} gives the proofs of Lemma~\ref{squarepairlemma1} and~\ref{strongsqpairlemma}.  Finally, in Section~\ref{sec:conc}, we give some concluding remarks and a summary of known results.

%For this paper, we focus on improving the lower bound for the minimum density of
%a 2-identifying code in both the hexagonal and square grid.  It is possible to also
%use this technique to show that the density of a 3-identifying code is at least 3/25, %but the proof is long and the improvement is minor so we will exclude it here.  The %proof of this fact will appear in ~\cite{StantonPending}.  To improve further on the %square grid case, we use the discharging method.

\section{Definitions and General Lemmas}\label{generalsection}

Let $G_S$ denote the square grid.  Then $G_S$ has vertex set $V(G_S)=\Z\times\Z$ and
$$ E(G_S)=\{\{u,v\}: u-v\in \{(0, \pm 1),(\pm 1,0)\}\} , $$
where subtraction is performed coordinatewise.

Let $G_H$ represent the hex grid.  We will use the so-called ``brick wall'' representation, whence $V(G_H)=\Z\times\Z$ and
$$ E(G_H)=\{\{u=(i,j),v\}: u-v\in \{(0,(-1)^{i+j+1}),(\pm 1,0)\}\} . $$

Consider an $r$-identifying code $C$ for a graph $G=(V,E)$.  Let
$c,c'\in C$ be distinct.  If $I_r(v)=\{c,c'\}$ for some $v\in V(G)$
we say that \begin{enumerate}
\item $c'$ \textdef{forms a pair} (with $c$) and
\item $v$ \textdef{witnesses a pair} (that contains $c$).
\end{enumerate}

For $c\in C$, we define the set of witnesses of pairs that contain $c$.  Namely,
$$ P(c)=\{v:I_r(v)=\{c,c'\}, \text{ for some }c\neq c'\} . $$
We also define $p(c)=|P(c)|$.  In other words,
$P(c)$ is the set of all vertices that witness a pair containing $c$ and $p(c)$ is the number of vertices that witness a pair containing $c$. Furthermore, we call $c$ a \textdef{$k$-pair codeword} if $p(c)=k$.

We start by noting two facts about pairs which are true for any code on any graph.

\begin{fact}\label{subsetfact1}
    Let $c$ be a codeword and $S$ be a set of vertices such that for
    each $s\in S$, $s$ witnesses a pair containing $c$.  If $v\not\in S$ and
    $B_2(v)\subset\bigcup_{s\in S}B_2(s)$, then $v$ does not witness a pair containing $c$.
\end{fact}

\begin{proof}
Suppose $v$ witnesses a pair containing $c$.  Hence, $I_2(v)=\{c,c'\}$ for some $c'\neq c$.  Then $c'\in B_2(v)$ and so $c'\in B_2(s)$ for some $s\in S$.  But then $\{c,c'\}\subset I_2(s)$.  But since $I_2(s)\neq I_2(v)$, $|I_2(s)|>2$, contradicting the fact that $s$ witnesses a pair.  Hence $v$ does not witness such a pair.
\end{proof}

\begin{fact}\label{subsetfact2}
    Let $c$ be a codeword and $S$ be any set with $|S|=k$.  If $v\in S$ and $$B_2(v)\subset\bigcup_{s\in S\atop s\neq v}B_2(s)$$ then at most $k-1$ vertices in $S$ witness pairs containing $c$.
\end{fact}

\begin{proof}
The result follows immediately from Fact~\ref{subsetfact1}.  If each vertex in $S-\{v\}$ witnesses a pair, then $v$ cannot witness a pair.  Hence, either $v$ does not witness a pair or some vertex in $S$ does not witness a pair.
\end{proof}

Lemma~\ref{generalpairlemma} is a general statement about vertex-identifying codes and has a similar proof to Theorem 2 in~\cite{Karpovsky1998}.
\begin{lemma}\label{generalpairlemma}
   Let $C$ be an $r$-identifying code for the square or hex grid.  Let $p(c)\le k$ for any codeword.  Let $D(C)$ represent the density of $C$, then if $b_r=|B_r(v)|$ is the size of a ball of radius $r$ centered at any vertex $v$,
   $$ D(C)\ge\frac{6}{2b_r+4+k} . $$
\end{lemma}

\begin{proof}
   We first introduce an auxiliary graph $\Gamma$.  The vertices of $\Gamma$ are the vertices in $C$ and $c\sim c'$ if and only if $c$ forms a pair with $c'$.  Then we clearly have $\deg_\Gamma(c)=p(c)$.  Let $\Gamma[C\cap G_m]$ denote the induced subgraph of $\Gamma$ on $C\cap G_m$.  It is clear that if $\deg_\Gamma(c)\le k$ then $\deg_{\Gamma[C\cap G_m]}\le k$.

   The total number of edges in $\Gamma[C\cap G_m]$ by the handshaking lemma is $$ \frac12\sum_{c\in \Gamma[C\cap G_m]}\deg_{\Gamma[C\cap G_m]} \le (k/2)|C\cap G_m| .$$
   But by our observation above, we note that the total number of pairs in $C\cap G_m$ is equal to the number of edges in $\Gamma[C\cap G_m]$.  Denote this quantity by $P_m$.  Then
   $$ P_m\le (k/2)|C\cap G_m| . $$

   Next we turn our attention to the grid in question.  The arguments work for either the square or hex grid.  Note that if $C$ is a code on the grid, $C\cap G_m$ may not be a valid code for $G_m$.  Hence, it is important to proceed carefully.  Fix $m>r$.  By definition, $G_{m-r}$ is a subgraph of $G_m$.  Further, for each vertex $v\in V(G_{m-r})$, $B_r(v)\subset V(G_m)$.  Hence $C\cap G_m$ must be able to distinguish between each vertex in $G_{m-r}$.

   Let $n=|G_m|$ and $K=|C\cap G_m|$.  Let $v_1,v_2,v_3,\ldots, v_n$ be the vertices of $G_m$ and let $c_1,c_2,\ldots, c_K$ be our codewords.  We consider the $n\times K$ binary matrix $\{a_{ij}\}$ where $a_{ij}=1$ if $c_j\in I_r(v_i)$ and $a_{ij}=0$ otherwise.  We count the number of non-zero elements in two ways.

   On the one hand, each column can contain at most $b_r$ ones since each codeword occurs in $B_r(v_i)$ for at most $b_r$ vertices.  Thus, the total number of ones is at most $b_r\cdot K$.

   Counting ones in the other direction, we will only count the number of ones in rows corresponding to vertices in $G_{m-r}$.  There can be at most $K$ of these rows that contain a single one and at most $P_m$ of these rows which contain 2 ones.  Then there are $|G_{m-k}|-K-P_m$ left corresponding to vertices in $G_{m-k}$ and so there must be at least 3 ones in each of these rows.  Thus the total number of ones counted this way is at least $K +
   2P_m+3(|G_{m-r}|-K-P_m) = -2K +3|G_{m-r}|-P_m$.  Thus
   \begin{equation}   \label{eqn:mainineq} b_rK\ge-2K +3|G_{m-r}|-P_m .
  \end{equation}

   But since $P_m\le (k/2)K$, this gives
   $$ b_rK\ge-2K +3|G_{m-r}|-(k/2)K . $$
   Rearranging the inequality and replacing $K$ with $|C\cap G_m|$ gives $$ \frac{|C\cap G_m|}{|G_{m-r}|}\ge \frac{6}{2b_r +4+k} .$$

   Then
   \begin{eqnarray*}
      D(C) & = & \limsup_{m\rightarrow\infty}\frac{|C\cap G_m|}{|G_m|} \\
      & = & \limsup_{m\rightarrow\infty}\frac{|C\cap G_{m}|}{|G_{m-r}|}\cdot \limsup_{m\rightarrow\infty}\frac{|G_{m-r}|}{|G_m|} \\
      & \ge & \frac{6}{2b_r +4+k}\cdot \limsup_{m\rightarrow\infty}\frac{\frac32(m-r)^2 + \frac32 (m-r)+1}{\frac32m^2+\frac32m+1} \\
      & = & \frac{6}{2b_r +4+k} .
   \end{eqnarray*}
\end{proof}

\section{Lower Bound for the Hexagonal Grid}\label{hexsection}

Lemma~\ref{hex2lemma} establishes an upper bound of 6 for the degree of the graph $\Gamma$ formed by a code in the hex grid, which allows us to prove Theorem~\ref{hexr2theorem}.

\begin{lemma}\label{hex2lemma}
Let $C$ be a 2-identifying code for the hex grid. For each $c\in C$, $p(c)\le 6$.
\end{lemma}

\begin{proof}
    Let $C$ be a code and $c\in C$ is an arbitrary member.  Let $u_1,u_2,$ and $u_3$ be the neighbors of $c$ and let
    $\{u_{i1},u_{i2}\}=B_1(u_i)-\{u_i,c\}$.

    \noindent\textbf{Case 1:}  $|I_2(c)|\ge 2$

    There exists some $c'\in C\cap B_2(c)$ with $c'\neq c$.  Without
    loss of generality, assume that $c'\in \{u_1,u_{11},u_{12}\}$.
    Since $I_2(c),I_2(u_1),I_2(u_{11}),I_2(u_{12})\supseteq \{c,c'\}$
    at most one of $c,u_1,u_{11},u_{12}$ witnesses a pair containing
    $c$.

    Now, $p(c)\le 6$ unless each of $u_2,u_3,u_{21},u_{22},
    u_{31},u_{32}$ witnesses a pair.

    If $u_2$ and $u_3$ each witness a pair, then we have
    $u_i\not\in C$ for $i=1,2,3$; otherwise
    $I_2(u_2)=\{c,u_i\}=I_2(u_3)$ and so $u_2$ and $u_3$ are not
    distinguishable by our code.  Thus, there must be some $c''\in
    C\cap (B_2(u_2)-\{c,u_1,u_2,u_3\})$. This forces
    $c''\in  B_2(u_{21})\cup B_2(u_{22})$ and so either
    $\{c,c''\}\subseteq I_2(u_{21})$ or $\{c,c''\}\subseteq I_2(u_{22})$. Hence, one of these cannot witness a pair and still be
    distinguishable from $u_2$.  This ends case 1.

    \noindent\textbf{Case 2:}  $I_2(c) = \{c\}$

    First note that $c$ itself does not witness a pair.

    If $u_1$ witnesses a pair, then there is some $c''\in C\cap(
    B_2(u_1)-B_2(c)) \subseteq C\cap ( B_2(u_{11})\cup
    B_2(u_{12}))$ and so either
    $\{c,c''\}\subseteq I_2(u_{11})$ or $\{c,c''\}\subseteq I_2(u_{12})$
    and so one of these cannot witness a pair and still be
    distinguishable from $u_1$.  Hence at most two of
    $\{u_1,u_{11},u_{12}\}$ can witness a pair.

    Likewise at most at most two of
    $\{u_2,u_{21},u_{22}\}$ and
    $\{u_3,u_{31},u_{32}\}$ can witness a pair. Thus $p(c)\le 6$.
    This ends both case 2 and the proof of the lemma.
\end{proof}

\begin{proofcite}{Theorem~\ref{hexr2theorem}}
   Using Lemmas~\ref{generalpairlemma} and~\ref{hex2lemma}, if $C$ is a $2$-identifying code in the hexagonal grid, then
   $$ D(C)\geq\frac{6}{2b_2+4+6}=\frac{6}{30}=\frac{1}{5} . $$
\end{proofcite}

\section{Lower Bounds for the Square Grid}\label{squaresection}

Lemma~\ref{squarepairlemma1} establishes an upper bound of 8 for the degree of the graph $\Gamma$ formed by a code in the square grid, which allows us to prove Theorem~\ref{squarer2theorem}.  Then we prove Lemma~\ref{strongsqpairlemma}, which bounds the average degree of $\Gamma$ by 7, allowing for the improvement in Theorem~\ref{theorem:mainsquare}.

\begin{lemma}\label{squarepairlemma1}
Let $C$ be a 2-identifying code for the square grid. For each $c\in C$, $p(c)\le 8$.
\end{lemma}

\begin{proof}
Let $c\in C$, a $2$-identifying code in the square grid. Without loss of generality, we will assume that $c=(0,0)$.

\begin{figure}[ht]
    \centering
    \includegraphics{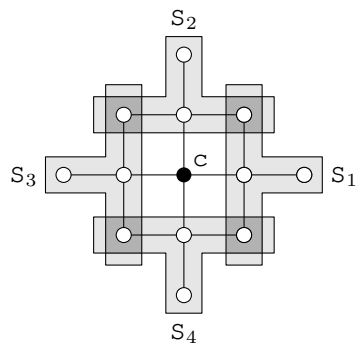}
    \caption{The sets $S_1$, $S_2$, $S_3$ and $S_4$.}\label{fig:ssets}
\end{figure}

\begin{figure}[ht]
\centering
\includegraphics{squarecenter.mps}
\caption {The ball of radius 2 around $c$.  A configuration of 9
vertices witnessing pairs is not possible if $|I_2(c)|= 2$.\newline $\bullet$~At most 7 of the vertices in gray triangles may
witness a pair.\newline $\bullet$~At most one of the vertices in white
triangles may witness a pair.}\label{fig:squarecenter}
\end{figure}

\textbf{Case 1:} $c$ witnesses a pair.

This case implies immediately that $|I_2(c)|=2$.  The other codeword in $I_2(c)$, namely $c'$, is in one of the
following 4 sets, the union of which is $B_2(c)-\{c\}$.  See Figure~\ref{fig:ssets}.
$$ \begin{array}{rrrrr}
      S_1:=\{ & (1,0),  & (1,1),  & (1,-1),  & (2,0)\} \\
      S_2:=\{ & (0,1),  & (1,1),  & (-1,1),  & (0,2)\} \\
      S_3:=\{ & (-1,0), & (-1,1), & (-1,-1), & (-2,0)\} \\
      S_4:=\{ & (0,-1), & (1,-1), & (-1,-1), & (0,-2)\}
\end{array}$$

If, however, $c'\in S_i$, then no $s\in S_i$ can witness a pair because $\{c,c'\}\subseteq I_2(s)$ and $s$ could not be distinguished from $c$. Without loss of generality, assume that $c'\in S_3$.  Thus, all vertices witnessing pairs in $I_2(c)$ are in the set
$$ R:=\left\{(x,y) : (x,y)\in B_2(c), x\geq 0\right\} . $$
But because
$$ B_2\left((1,0)\right)\subseteq\bigcup_{s\in S_1\cup\{c\}}B_2(s) , $$
Fact~\ref{subsetfact1} gives that not all members of $S_1\cup\{c\}$ can witness a pair. See Figure~\ref{fig:squarecenter}.

Therefore, $p(c)\leq 8$ and, without loss of generality, $c'\in S_3$ and at least one of $S_1$ does not witness a pair.  This ends Case 1.

\textbf{Case 2:} $c$ does not witness a pair.

This case implies immediately that either $|I_2(c)|\geq 3$ or $I_2(c)=\{c\}$.

First suppose $|I_2(c)|\ge 3$. There must be two distinct codewords $c',c''\in S_1\cup S_2 \cup S_3 \cup S_4$. If $c',c''$ are in the same set $S_i$ for some $i$, then $\{c,c',c''\}\subset I_2(s)$ for any $s\in S_i$
and so no vertex in $S_i$ witnesses a pair. Thus, the only vertices which can witness a pair are in $B_2(c)-(S_i\cup \{c\})$.  There are only 7 of these, so $p(c)\le 7$. (See the gray vertices in Figure~\ref{fig:squarecenter}).

If $c'\in S_i$ and $c''\in S_j$ for some $i\neq j$, then only one vertex in each of $S_i$ and $S_j$ can witness a pair.  There are at most 5 other vertices not in $S_i\cup S_j-\{c\}$ and so $p(c)\le 7$.

Thus, if $|I_2(c)|\ge 3$, then $p(c)\le 7$.

\begin{figure}[ht]
\centering
\includegraphics{squaretriangle1.mps}
\caption{A right angle of witnesses.  \newline $\bullet$~Black circles indicate codewords.
\newline $\bullet$~White circles indicate non-codewords.
\newline $\bullet$~Gray triangles indicate vertices that witness a pair.  \newline $\bullet$~White triangles indicate vertices that do not witness a pair.  \newline No vertices in $B_2(c)-\{c\}$ can be codewords, neither can those which are distance no more than 2 from two vertices in this right angle of witnesses.}\label{fig:squaretriangle1}
\end{figure}

\begin{figure}[ht]
\centering
\includegraphics{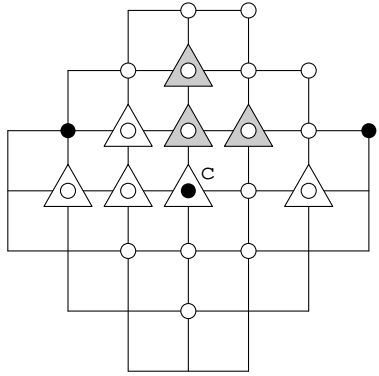}
\caption{A right angle of witnesses, continuing from Figure~\ref{fig:squaretriangle1}.  Let $c=(0,0)$. Vertices $(-2,1)$ and $(3,1)$ must be codewords and so none of  \nobreak{$\{(-1,1),(-1,0),(-2,0),(2,0)\}$}
can witness pairs.}\label{fig:squaretriangle2}
\end{figure}

Second, suppose $I_2(c)=\{c\}$.  We will define a \textdef{right angle of witnesses} to be subsets of 3 vertices of $I_2(c)$ that all witness pairs and are one of the following 8 sets: $\{(1,0),(2,0),(1,\pm 1)\}$, $\{(0,1),(0,2),(\pm 1,1)\}$, $\{(-1,0),(-2,0),(-1,\pm 1)\}$, and \linebreak $\{(0,-1),(0,-2),(\pm 1,-1)\}$.
If a right angle is present then, without loss of generality, let it be $\{(0,1),(0,2),(1,0)\}$.  See Figure~\ref{fig:squaretriangle1}.  In order for these all to be witnesses, then $I_2((0,1))$ must have one codeword not in $B_2((0,2))\cup B_2((1,1))$, which can only be $(-2,1)$.  Since $\{(0,0),(-2,1)\}\subseteq B_2((-1,1)),B_2((-1,0)),B_2((-2,0))$, none of those three vertices can witness a pair.

In addition, $I_2((1,1))$ must contain a codeword not in $B_2((0,1))\cup B_2((0,2))$, which can only be $(3,1)$.  See Figure~\ref{fig:squaretriangle2}. Since $\{(0,0),(3,1)\}\subseteq B_2((2,0))$, the vertex $(2,0)$ cannot witness a pair.

Finally, it is not possible for all of $(-1,-1),(0,-1),(1,-1),(0,-2)$ to be witnesses because the only member of $B_2((0,-1))$ that is not in the union of the second neighborhoods of the others is the vertex $(0,1)$, which cannot be a codeword in this case.  Hence, at most 7 members of $B_2(c)$ can witness a pair if $B_2(c)$ has a right angle of witnesses.

Consequently, if $c$ does not witness a pair and $p(c)\geq 8$, then $I_2(c)=\{c\}$ and $B_2(c)$ fails to have a right angle of witnesses.  We can enumerate the remaining possibilities according to how many of the vertices $\{(1,1),(-1,1),(-1,-1),(1,-1)\}$ are witnesses.  If 1, 2 or 3 of them are witnesses and there is no right angle of witnesses, it is easy to see that there are at most 7 witnesses in $B_2(c)$ and so $p(c)\leq 7$.

The first remaining case is if 0 of them are witnesses, implying each of the eight vertices $(\pm 1,0)$, $(\pm 2,0)$, $(0,\pm 1)$ and $(0,\pm 2)$ are witnesses.  The second remaining case is if 4 of them are witnesses.  This implies that at most one of $\{(1,0),(2,0)\}$ are witnesses and similarly for $\{(0,1),(0,2)\}$, $\{(-1,0),(-2,0)\}$ and $\{(0,-1),(0,-2)\}$.

This ends both Case 2 and the proof of the lemma.  So, $p(c)\leq 8$ with equality only if one of two cases in the previous paragraph holds.
\end{proof}

\begin{proofcite}{Theorem~\ref{squarer2theorem}}
   Using Lemmas~\ref{generalpairlemma} and~\ref{squarepairlemma1}, if $C$ is a $2$-identifying code in the square grid, then
   $$ D(C)\geq\frac{6}{2b_2+4+8}=\frac{6}{38}=\frac{3}{19} . $$
\end{proofcite}

\begin{lemma}\label{strongsqpairlemma}
  Let $G_m$ denote the induced subgraph of the square grid on the
  vertex set $[-m,m]\times[-m,m]$.  Let $C$ be a code for the square grid.
  Then $\sum_{c\in C\cap G_m}p(c)\le 7|C\cap G_m|$.
\end{lemma}

\begin{proof}
  Let $C$ be a code on $G_S$ and let
  $$ R(c)=\{c': I_2(v)=\{c,c'\}\text{ for some } c\in C\}. $$
  Suppose that $p(c)=8$ for some $c\in C$.  We claim that one of the two following properties holds.
  \begin{enumerate}
    \item[(P1)] There exist distinct $c_1,c_2,c_3\in R(c)$
    such that $p(c_1)\le 4$ and $p(c_i) \le 6$ for $i=2,3$.
    \item[(P2)] There exist distinct $c_1,c_2,c_3,c_4,c_5,c_6\in R(c)$
    such that $p(c_i)\le 6$ for all $i$.
  \end{enumerate}

  We will prove this by characterizing all possible $8$-pair
  vertices, but first we wish to define 3 different types of
  codewords.  The definition of each type extends by taking translations and rotations. So, we may assume in defining the types that $c=(0,0)$.

  We say that $c$ is a \textdef{type 1} codeword if $(0,1),(0,-1)\in C$. See Figure~\ref{fig:vertextype1}.

  We say that $c$ is a \textdef{type 2} codeword if $(-1,2),(2,-1)\in C$. See Figure~\ref{fig:vertextype2}.

  We say that $c$ is a \textdef{type 3} codeword if $(-2,1),(2,1)\in C$. See Figure~\ref{fig:vertextype3}.

Claim~\ref{squareclaim1} shows that adjacent codewords do not need to be considered because they are in few pairs.

\begin{claim}\label{squareclaim1}  If $c$ is adjacent to another codeword, then $p(c)\le 6$.\end{claim}

\begin{proof} Without loss of generality, assume that $c=(0,0)$ and that $(0,1)$ is a codeword.  Then
$$ (-1,0), (0,0), (0,1), (0,2), (1,0), (1,1), (-1,0), (-1,1)$$
are all at most distance 2 from both codewords and so at most 1 of them can witness a pair.  Thus, the other 7 do not witness pairs containing $c$.  Since $|B_2(c)|=13$, $p(c)\le 13-7=6$.  This proves Claim~\ref{squareclaim1}.
\end{proof}

Claims~\ref{claimtype1},~\ref{claimtype2} and~\ref{claimtype3} show that types 1, 2 and 3 codewords, respectively, are not in many pairs.

\begin{claim}\label{claimtype1}
  If $c$ is a type 1 codeword, then $p(c)\le 4$.
\end{claim}

\begin{figure}[ht] \centering
\includegraphics{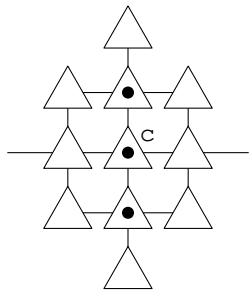}
\caption{Vertex $c$ is a type 1 codeword.  At most 2 of the 11 vertices marked by triangles can witness a pair.}\label{fig:vertextype1}
\end{figure}

\begin{proof} Without loss of generality, let $c=(0,0)$.  We consider all vertices which are distance 2 from $c$ and either $(0,1)$ or $(0,-1)$.  There are 11 such vertices and at most 2 of them can witness pairs, so $p(c)\le 4$.  See
Figure~\ref{fig:vertextype1}. This proves Claim~\ref{claimtype1}.
\end{proof}

\begin{claim}\label{claimtype2}
   If $c$ is a type 2 codeword, then $p(c)\le 6$.
\end{claim}

\begin{figure}[ht]
    \centering
    \includegraphics{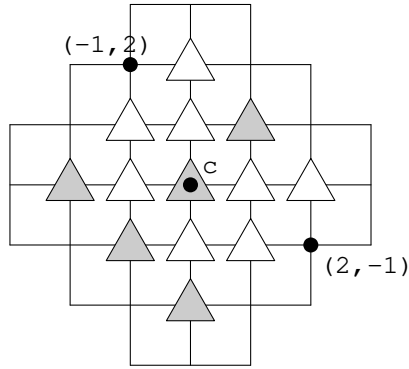}
    \caption{Vertex $c=(0,0)$ is a type 2 codeword.  At most 2 of the 8 vertices marked by white triangles can witness pairs.  At most 4 of the 5 vertices marked by gray triangles can witness pairs.}\label{fig:vertextype2}
\end{figure}

\begin{proof} Without loss of generality, let $c=(0,0)$.  We consider all vertices which are distance at most 2 from $c$ and distance at most 2 from either $(-1,2)$ or $(2,-1)$.  There are 8 such vertices and at most 2 of them can witness pairs.  The remaining 5 vertices are $c$ and the vertices in the set $S=\{(-2,0),(-1,-1),(0,-2),(1,1)\}$.  But then $B_2(c)\subset\bigcup_{s\in S}B_2(s)$ and, by Fact~\ref{subsetfact1} at most 4 of those remaining 5 vertices can witness pairs. Thus, $p(c)\le 6$.  See Figure~\ref{fig:vertextype2}. This proves Claim~\ref{claimtype2}.
\end{proof}

\begin{claim}\label{claimtype3}
   If $c$ is a type 3 codeword, then $p(c)\le 6$.
\end{claim}

\begin{figure}[ht]
    \centering
    \includegraphics{vertextype3.mps}
    \caption{Vertex $c=(0,0)$ is a type 3 codeword. \newline $\bullet$~$T_0$ vertices are black. \newline $\bullet$~$T_1$ vertices are
    white. \newline
    $\bullet$~$T_2$ vertices are marked by diagonal lines. \newline $\bullet$~$T_3$ vertices are gray.}\label{fig:vertextype3}
\end{figure}

\begin{proof} Without loss of generality, let $c=(0,0)$. We partition $B_2(c)- \{c\}$ into 4 sets:
$$ \begin{array}{rrrrr}
      T_0:= \{ &          &         & (0,1),  & (0,2)\} \\
      T_1:= \{ &          & (-2,0), & (-1,0), & (-1,1)\} \\
      T_2:= \{ &          &  (2,0), &  (1,0), &  (1,1)\} \\
      T_3:= \{ & (-1,-1), & (0,-1), & (1,-1), & (0,-2)\} \\
\end{array}$$

At most 1 vertex in $T_0$ witnesses a pair since $|I_2(0,1)|\ge 3$.

At most 1 vertex in $T_1$ can witness a pair since every vertex in $T_1$ is at most distance 2 from $(-2,1)$. Likewise, at most 1 vertex in $T_2$ can witness a pair.

If all vertices in $T_3$ witness pairs, then $I_2((0,-1)) = \{(0,0),(0,1)\}$ since $(0,1)$ is the only vertex in $B_2((0,-1))$ which is not in $B_2(s)$ for any other $s\in T_3$.  But then $c$ is adjacent to another codeword, and by Claim~\ref{squareclaim1}, $p(c)\le 6$.  So we may assume that at most 3 vertices in $T_3$ form pairs with $c$.

Now, if $c$ does not itself witness a pair, these partitions give $p(c)\le 6$.  If $c$ does witness a pair, then there must be another codeword $c'\in S_i$ for some $i$.  But then we see that no other vertex in $S_i$ can witness a pair, since every vertex in $S_i$ is at most distance two from  $c'$.  Thus, $p(c)\le 6$.  See    Figure~\ref{fig:vertextype3}.  This proves Claim~\ref{claimtype3}.
\end{proof}

We are now ready to characterize the 8-pair codewords.  

\begin{claim}\label{claimpair}
   If $c\in C$ witnesses a pair and $p(c)=8$, then $c$ satisfies property (P1).
\end{claim}

\begin{figure}[ht]
   \centering
   \includegraphics{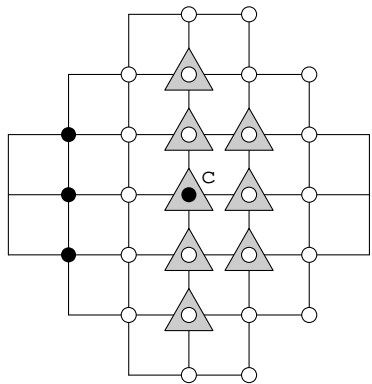}
   \caption{Codeword $c=(0,0)$ witnesses a pair and is an $8$-pair codeword.  The gray triangles are vertices that form pairs with $c$.  Vertex $(-2,0)$ is a type 1 codeword.}\label{fig:squarecpair}
\end{figure}

\begin{proof} Without loss of generality, let $c=(0,0)$.  
Recall Case 1 of the proof of Lemma~\ref{squarepairlemma1}.  That is, $p(c)\leq 8$ and, without loss of generality, equality implies that there is a $c'\in C\cap S_3$ and at least one of $S_1=\left\{(1,-1),(1,0),(1,1),(2,0)\right\}$ does not witness a pair.

If $p(c)\leq 7$, the proof is finished, so let us assume that $p(c)=8$ and hence exactly one of the vertices in $S_1$ does not witness a pair.  We will show that it is $(2,0)$.  So, suppose that $(1,y)$ does not witness a pair. Recall that $R=\left\{(x,y) : (x,y)\in B_2(c), x\geq 0\right\}$.

If $y\in\{-1,1\}$, then
$$ B_2\left((1,0)\right)\subseteq\bigcup_{s\in R-\{(1,y),(1,0)\}}B_2(s) $$
and, by Fact~\ref{subsetfact1}, neither $(1,y)$ nor $(1,0)$ witnesses a pair and $p(c)\leq 7$.

If $y=0$, then
$$ B_2\left((1,1)\right)\subseteq\bigcup_{s\in R-\{(1,0),(1,1)\}}B_2(s) $$
and, by Fact~\ref{subsetfact1}, neither $(1,0)$ nor $(1,1)$ witnesses a pair and $p(c)\leq 7$. It follows that each vertex in $R'=R-\{(2,0)\}$ must witness a pair containing $c$.

Each vertex which is distance 2 or less from 2 vertices in $R'$ cannot be a codeword.  Thus, $(-2,0)$ is the only vertex in $B_2(c)$ other than $c$ which has not been marked as a non-codeword and so $(-2,0)\in C$.  Since $(0,0)\in C$, the vertex $(-2,1)$ is the only possibility for a second codeword for $(0,1)$ and $(-2,-1)$ is the only possibility for a second codeword for $(0,-1)$.  See Figure~\ref{fig:squarecpair}.

Then $(-2,0)$ is a type 1 codeword and so it is in at most 4 pairs.  Codewords $(-2,1)$ and $(-2,-1)$ are both adjacent to another codeword, so they are in at most 6 pairs.  Hence, $c$ satisfies Property (P1).  This proves Claim~\ref{claimpair}.
\end{proof}

\begin{claim}\label{claimnopair}
   If $c\in C$ does not witness a pair and $p(c)=8$, then $c$ satisfies either property (P1) or property (P2).
\end{claim}

\begin{proof} Without loss of generality, let $c=(0,0)$.
Recall Case 2 of the proof of Lemma~\ref{squarepairlemma1}.  
That is, $p(c)\leq 8$ and, without loss of generality, equality implies $I_2(c)=\{c\}$.  Furthermore, one of the following two cases occurs: \\
(1) The eight witnesses are the vertices $(\pm 1,0)$, $(\pm 2,0)$, $(0,\pm 1)$ and $(0,\pm 2)$. \\
(2) The witnesses include $\{(1,1),(-1,1),(-1,-1),(1,-1)\}$ as well as exactly one of each of the following pairs: $\{(1,0),(2,0)\}$, $\{(0,1),(0,2)\}$, $\{(-1,0),(-2,0)\}$ and \\
$\{(0,-1),(0,-2)\}$.

If case (1) occurs, then the eight witnesses are the vertices $(\pm 1,0)$, $(\pm 2,0)$, $(0,\pm 1)$ and $(0,\pm 2)$.  In this case, simply observe that $B_2((1,0))$ is a subset of the other seven witnesses.  This contradicts Fact~\ref{subsetfact2} and so this case cannot occur.

\begin{figure}[ht]
    \centering
    \includegraphics{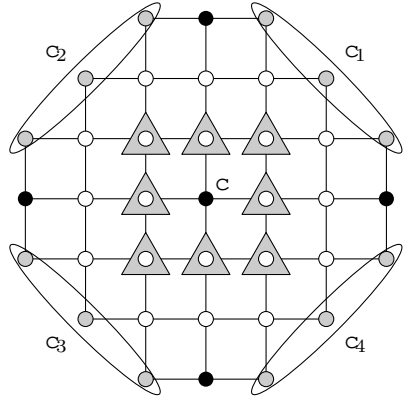}
    \caption{Codeword $c=(0,0)$ fails to witness a pair and is an $8$-pair codeword.  Exactly one of the gray vertices in each oval is a codeword.}\label{fig:square8cw}
\end{figure}

So, we may assume that case (2) occurs.  The vertex $(2,1)$ cannot be a codeword because $\{(0,0),(2,1)\}\subseteq B_2((1,1)),B_2((1,0)),B_2((2,0))$ and so at most one of these three vertices witness pairs, a contradiction to case (2).  By symmetry, none of the following vertices witness pairs: $$ (2,1),(1,2),(-1,2),(-2,1),(-2,-1),(-1,-2),(1,-2),(2,-1) . $$

In order to distinguish $(1,0)$ from $(0,0)$, the only vertex available to be a codeword is $s_1:=(3,0)$ and symmetrically, $s_2:=(0,3)$, $s_3:=(-3,0)$ and $s_4:=(0,-3)$ are codewords.  This implies that each of $(1,0)$, $(0,1)$, $(-1,0)$ and $(0,-1)$ witness pairs.

Then, for the other 4 pairs, there are exactly 3 choices for codewords which are not in the ball of radius 2 for any of our other pairs. See Figure~\ref{fig:square8cw}.
$$ \begin{array}{c|rcl}
      \text{Vertex} & \multicolumn{3}{c}{\text{Other Codeword}} \\
      \hline
      (1,1) & c_1 & \in & \{(3,1),(2,2),(1,3)\} \\
      (-1,1) & c_2 & \in & \{(-1,3),(-2,2),(-3,1)\}  \\
      (-1,-1) & c_3 & \in & \{(-3,-1),(-2,-2),(-1,-3)\}  \\
      (1,-1) & c_4 & \in & \{(1,-3),(2,-2),(3,-1)\}
    \end{array} $$

For each $c_i$, either $c_i$ is adjacent to another codeword or $c_i$ is a type 2 codeword.  Claims~\ref{squareclaim1} and~\ref{claimtype2} imply that, in either case, $p(c_i)\le 6$.  It remains to show that one of the following holds:  (1) There exist $i\neq j$ such that $p(s_i)\le 6$ and $p(s_j)\le 6$, hence $c$ satisfies (P1). (2) There exists an $i$ such that $p(s_i)\le 4$, hence $c$ satisfies (P2).

First, suppose that there are $c_i$, $c_j$, $i\neq j$ such that $c_i\sim s_k$, $c_j\sim s_\ell$.  If $k=\ell$, then $s_k$ is a type 1 codeword and so $p(s_k)\le 4$.  If $k\neq \ell$, then both $s_k$ and $s_\ell$ are adjacent to another codeword and so $p(s_k)\le 6$ and $p(s_\ell)\le 6$.  Either (P1) or (P2) is satisfied, respectively.

If there is at most one $c_i$ such that $c_i\sim s_k$ for some $k$, then we have three codewords of the form $(\pm 2,\pm 2)$.  Without loss of generality, assume that $(2,2),(2,-2),$ and $(-2,2)$ are codewords.  In this case, $(3,0)$ and $(0,3)$ are type 3 codewords and hence $p((3,0))\le 6$ and $p((0,3))\le 6$.  So again, (P2) is satisfied.

This proves Claim~\ref{claimnopair}.
\end{proof}

Finally, we can finish the proof of Lemma~\ref{strongsqpairlemma} by way of the discharging method. (For a more extensive application of the discharging method on vertex identifying codes, see Cranston and Yu~\cite{Cranston2009}.)  Let $\Gamma$ denote an auxiliary graph with vertex set $C\cap G_m$ for some $m$. There is an edge between two vertices $c$ and $c'$ if and only if $I_2(v)=\{c,c'\}$ for some $v\in V(G_S)$. For each vertex $v$ in our auxiliary graph $\Gamma$, we assign it an initial charge of $d(v)-7$.  Note that $\sum_{c\in C\cap G_m} p(c)-7 = \sum_{v\in \Gamma}\deg_{\Gamma}(v)-7$.  We apply the following discharging rules if $\deg_{\Gamma}(v)=8$.
\begin{enumerate}
  \item If $v$ is adjacent to one vertex of degree at most 4 and two of degree at most 6 (condition (P1)), then discharge 2/3 to a vertex of degree at most 4 and 1/6 to two vertices of degree at most 6.
  \item If $v$ is adjacent to 6 vertices of degree at most 6 (condition (P2)), then discharge 1/6 to 6 neighbors of degree at most 6.
\end{enumerate}

We have proven that one of the above cases is possible. Let $e(v)$ be the charge of each vertex after discharging takes place.  We show that $e(v)\le 0$ for each vertex in $\Gamma$.

If $\deg_{\Gamma}(v)=8$, then our initial charge was $1$.  In either of the two cases, we are discharging a total of 1 unit to its neighbors.    Since no degree 8 vertex receives a charge from any other vertex, we have $e(v)=0$.

If $d(v)=7$ then its initial charge is 0 and it neither gives nor receives a charge and so $e(v)=0$.

If $5\le \deg_{\Gamma}(v) \le 6$, then its initial charge was at most $-1$. Since this vertex has at most 6 neighbors and can receive a charge of at most $1/6$ from each of them, this gives $e(v)\le 0$.

If $\deg_{\Gamma}(v)\le 4$, then its initial charge was at most $-3$. Since this vertex has at most 2 neighbors and can receive a charge of at most $2/3$ from each of them, this gives $e(v)\le -1/3<0$.

Since no vertex can have degree more than 8, this covers all of the cases.  Then we have
$$ \sum_{c\in C\cap G_m} (p(c)-7)=\sum_{v\in \Gamma}\left(\deg_{\Gamma}(v)-7\right)=\sum_{v\in \Gamma}e(v) \le 0 . $$
Therefore, it follows that $\sum_{c\in C\cap G_m} p(c)\le \sum_{c\in C\cap G_m} 7=7|C\cap G_m|$.
\end{proof}

\begin{proofcite}{Theorem ~\ref{theorem:mainsquare}}
Consider $G_m$ and let $C$ be a code for $G_S$ and $C\cap G_m=\{c_1,c_2,\ldots, c_K\}$.  Recall inequality (\ref{eqn:mainineq}) from Theorem~\ref{generalpairlemma}.  In this case, $b_2=13$ and Lemma~\ref{strongsqpairlemma} shows that $$ P_m\le \frac12\sum_{c\in C\cap G_m}p(c)\le\frac 72|C\cap G_m| .$$

Substituting the above inequality into inequality~(\ref{eqn:mainineq}) and rearranging gives
$$ \frac{|C\cap G_m|}{|G_{m-r}|}\ge \frac {6}{37} . $$  Taking the limit as $m\rightarrow\infty$ gives the desired $D(C)\ge 6/37$, completing the proof.
\end{proofcite}

\section{Conclusions}
\label{sec:conc}

Below is a table noting our improvements.

$$\begin{array}{|c|c||c|c|}
  \hline
  \multicolumn{4}{|c|}{\text{Hex Grid}}\\
  \hline
  r & \text{previous lower bounds} & \text{new lower bounds} & \text{upper bounds} \\
  \hline
  2 & 2/11\approx 0.1818 ^\text{~\cite{Karpovsky1998}} & 1/5 = 0.2 & 4/19\approx 0.2105 ^\text{~\cite{Charon2002}}\\
  \hline
  3 & 2/17\approx 0.1176 ^\text{~\cite{Charon2001}} & 3/25 = 0.12^\text{~\cite{StantonPending}} & 1/6\approx 0.1667 ^\text{~\cite{Charon2002}}\\
  \hline
  \multicolumn{4}{|c|}{\text{Square Grid}}\\
  \hline
  2 & 3/20=0.15 ^\text{~\cite{Charon2001}} & 6/37\approx 0.1622 & 5/29\approx 0.1724 ^\text{~\cite{Honkala2002}}\\
  \hline
\end{array}$$

This technique works quite well for small values of $r$, but we note that $b_r=|B_r(v)|$ grows quadratically in $r$, so the denominator in Lemma~\ref{generalpairlemma} would grow quadratically.  But the known the lower bounds for $r$-identifying codes is proportional to $1/r$ in all of the well-studied grids (square, hexagonal, triangular and king).  Therefore, our technique is less effective as $r$ grows.  

\bibliographystyle{plain}
\bibliography{bibtex}
\end{document}